\documentclass[11pt]{article}
\usepackage{amssymb}
\usepackage{amsmath}
\usepackage{amsthm}
\usepackage{latexsym}
\usepackage{amsfonts}
\usepackage{graphicx}
\usepackage{graphics}
\usepackage[french,english]{babel}

\theoremstyle{plain}
\newtheorem{thm}{Theorem}[section]   
\newtheorem{prop}[thm]{Proposition}
\newtheorem{lem}[thm]{Lemma}

\theoremstyle{definition}
\newtheorem{rem}[thm]{Remark}

\newtheorem*{Proof}{Proof}

\begin{document}
\title{\bf UNIVERSAL PAD\'E APPROXIMANTS ON SIMPLY CONNECTED DOMAINS}
\author{N. Daras$^1$, G. Fournodavlos$^2$, V. Nestoridis$^3$}

\author{N. Daras\thanks{Mathematics Department, Hellenic Military Academy. Email: ndaras@sse.gr}\and
G. Fournodavlos\thanks{Mathematics Department, University of Toronto. Email: grifour@math.toronto.edu}\and
 V. Nestoridis\thanks{Mathematics Department, University of Athens. Email: vnestor@math.uoa.gr}}

\date{}
\maketitle

%
\noindent
%
\begin{abstract}
The theory of universal Taylor series can be extended to the
case of Pad\'e approximants where the universal approximation is
not realized by polynomials any more, but by rational functions,
namely the Pad\'e approximants of some power series.
We present the first generic result in this direction,
for Pad\'e approximants corresponding to Taylor developments
of holomorphic functions in simply connected domains.
The universal approximation is required only on compact sets ${K}$
which lie outside the domain of definition and have connected complement.
If the sets ${K}$ are additionally disjoint from the boundary of the domain of definition,
then the universal functions can be smooth on the boundary.
\end{abstract}
{\em Subject Classification MSC2010}\,: primary 41A21, 30K05
secondary 30B10, 30E10, 30K99, 41A10, 41A20. \vspace*{0.2cm} \\
{\em Key words}\,: Universality, Baire's theorem, Pad\'e approximant
\section{Introduction}\label{sec1}
\noindent

Henri Pad\'e, in his thesis, under Charles Hermite,
arranged his approximants in a double array, now known as the Pad\'e table
\[{\left({\left[{{ p}}/{{ q}}\right]}_{{ f},{ \zeta}}\left({ z}\right)\right)}_
{{p},{q}\in\{{0},{1},{2},\ldots\}},\]
of the formal power series$\ \ { f}{ =}$ $ \sum^{\infty }_{{ \nu}{ =}{ 0}}{{{ \alpha}}_{{ \nu}}}{\left({ z}{ -}{ \zeta}\right)}^{{ \nu}}=\left({{ \alpha}}_{{ 0}}{\rm ,\ }{{ \alpha}}_{{ 1}},{\rm \dots }\right)$. For$\ { q}{ =}{ 0}$, the Pad\'e approximants are polynomials and coincide with the partial sums of the power series. However, for$\ \ { q}\ge { 1}$ the Pad\'e approximants may exist or not and if they exist, they are rational functions.

Pad\'e approximants are used widely in pure and applied mathematics and they also have a wide range of applications in physics, chemistry and electronics; for more details, the reader may consult \cite{Bak1}, \cite{Dar1}, \cite{Mag}. In mathematics, the application which brought them to prominence in the 60's and 70's was regarding the localization of the singularities of functions in various problems \cite{Brui}. For example in inverse scattering theory, one would have a means to compute the coefficients of a power series $f$, which in turn could be used to find some Pad\'e approximant of $f$ and ultimately use its poles as predictors of the location of poles or other singularities of $f$. Under certain conditions on $f$, which were often satisfied in physical examples, this process could be theoretically justified.

The theory of universal Taylor series \cite{Chu,Luh,Nest4} can be extended to the case of Pad\'e approximants. In the context of universal Taylor series, in a simply connected domain ${ \Omega}$ \cite{Nest3,Nest1,Nest2} the partial sums of the Taylor development of a holomorphic function ${ f}\in {\mathcal H}\left({ \Omega}\right)$ approximate all polynomials uniformly on compact sets ${ K}$ disjoint from${ \ }{ \Omega}$ (or${ \ }\overline{{ \Omega}}\ $) with connected complement. Here, ${\mathcal H}\left({ \Omega}\right)$ denotes as usual the space of all functions holomorphic in $\Omega$, endowed with the topology of uniform convergence on compacta. If instead of partial sums we use a family of Pad\'e approximants${{ \ }\left[{{ p}}/{{ q}}\right]}_{{ f},{ \zeta}}\left({ z}\right)$ for the preceding approximation, where $\zeta$ belongs in $\Omega$ and ${ (}{ p},{ q}{ )}$ range in a subset $\mathfrak{F}\subset {{\mathbb N}}^{{ 2}}$ containing a sequence$\ \ \left({\widetilde{{ p}}}_{{ n}},{\widetilde{{ q}}}_{{ n}}\right)$,${ \ }{ n}{ =}{ 0},{ 1},{ 2}{ ,\dots }$ with$\ \ {\widetilde{{ p}}}_{{ n}}\to { +}\infty $, then we obtain the following generic result in ${\mathcal H}\left({ \Omega}\right)$ for any simply connected domain ${ \Omega}\subset {\mathbb C}$.
\begin{thm}\label{thm1.1}
 For the generic function${ \ \ }{ f}\in {\mathcal H}\left({ \Omega}\right)$, that is for every function in a dense and $G_\delta$ subset of $H(\Omega)$, subsequences of the Pad\'e approximation sequence of the Taylor series of f with center $\zeta$ in $\Omega$
\begin{align*} 
 \left({\left[{{{ p}}_{{ n}}}/{{{ q}}_{{ n}}}\right]}_{{ f},{ \zeta}}\left({ z}\right){ :\ \ }{ n}{ =}{ 0},{ 1},{ 2}{ ,\dots }\right),\ \ \left({{ p}}_{{ n}},{{ q}}_{{ n}}\right)\in \mathfrak{F}, 
 \end{align*}
 uniformly approximate any polynomial on any compact set $ { K}$ disjoint from 
$\Omega$, having connected complement.
\end{thm}
\begin{rem}\label{rem1.1}
We notice that in the previous theorem the centre ${ \zeta}\in { \Omega}$ may be fixed or the approximations can even be uniform when ${ \zeta}$ varies in compact subsets of ${ \Omega}$. Furthermore, if the compact sets ${ K}$ are disjoint from the boundary of$\ \ { \Omega}$, then the universal function${ \ \ }{ f}$ can be chosen to extend continuously on${ \ }\overline{{ \Omega}}$. This is a generic result in the space${ \ \ }{ A}\left({ \Omega}\right)$ of all holomorphic functions ${ f}\in {\mathcal H}\left({ \Omega}\right)$ having a continuous extension on the closure $\overline{{ \Omega}}$ of ${ \Omega}$ in ${\mathbb C}$, provided that $\left\{\infty \right\}\bigcup \left[{\mathbb C}\backslash \overline{{ \Omega}}\right]$ is connected.
\end{rem}
The above theorem
extends results (${q_n=0}$) that have been achieved, by several authors,
on universal Taylor series over the past years; see \cite{Nest2} and the references therein.

We note that a first description of our results was given in \cite{Dar2} with more restrictive conditions. In \cite{Four} we introduced a subtler method of proof, where we studied the problem of approximating the function $f$ itself by its Pad\'e approximants inside $\Omega$. Using this method, we were able to obtain refined generic results for universal Pad\'e approximants, in the form given in the present paper.

\section{PRELIMINARIES}\label{sec2}
\noindent

Let ${ \zeta}\in {\mathbb C}$ be fixed and ${ f}{ =}\displaystyle\sum^\infty_{{v=0}} {a_ v}\left({ z}{ -}{ \zeta}\right)^{v}$ be a formal power series (${{ a}}_{{ v}}{ =}{{ a}}_{{ v}}\left({ f},{ \zeta}\right)\in {\mathbb C}$). Often this series comes from the Taylor expansion of a holomorphic function on an open set containing$\ \ { \zeta}$. For$\ \ { p},{ q}\in \left\{{ 0},{ 1},{ 2}{ ,\dots }\right\}$, we define the Pad\'e approximant
\[{\left[{{ p}}/{{ q}}\right]}_{{ f},{ \zeta}}\left({ z}\right)\]
to be a rational function$\ \ \Phi\left({ z}\right){ =}{{\mathcal A}\left({ z}\right)}/{{\mathcal B}\left({ z}\right)}$, where ${\mathcal A}$ and ${\mathcal B}$ are polynomials satisfying$\ \ { deg}{\mathcal A}\le { p}$, ${ deg}{\mathcal B}\le { q}$, ${\mathcal B}\left({{ \zeta}}\right)= 1$, whose Taylor expansion \\
${ \ }\Phi\left({ z}\right){ =}\sum^{\infty }_{{ v}{ =}{ 0}}{{{ b}}_{{ v}}}{\left({ z}{ -}{ \zeta}\right)}^{{ v}}$ is such that
\begin{center}
${{ b}}_{{ v}}{ =}{{ a}}_{{ v}}$ for all ${ v}{ =}{ 0}{ ,\ }{ 1}{ ,\ }{ \dots }{ ,\ }{ p}{ +}{ q}$.
\end{center}
\noindent For${ \ }{ q}{ =}{ 0}$, the Pad\'e approximant ${\left[{{ p}}/{{ 0}}\right]}_{{ f},{ \zeta}}\left({ z}\right)$ exists trivially and it is equal to
\[{\left[{{ p}}/{{ 0}}\right]}_{{ f},{ \zeta}}\left({ z}\right){ =}\sum^{{ p}}_{{ v}{ =}{ 0}}{{{ a}}_{{ v}}}{\left({ z}{ -}{ \zeta}\right)}^{{ v}}.\]
For${ \ }{ q}\ge { 1}$, it is not always true that ${\left[{{ p}}/{{ q}}\right]}_{{ f},{ \zeta}}\left({ z}\right)$ exists, but if it exists it is unique as a rational function. A necessary and sufficient condition for existence and uniqueness of the polynomials ${\mathcal A}$ and ${\mathcal B}$ above is that the following ${ q}\times { q}$ Hankel determinant is different from zero:
\[{ det}\underbrace{\left( \begin{array}{c}
 \begin{array}{ccc}
{{ a}}_{{ p}{ -}{ q}{ +}{ 1}} &  \begin{array}{ccc}
{{ a}}_{{ p}{ -}{ q}{ +}{ 2}} & \cdots  \end{array}
 & {{ a}}_{{ p}}{ \ \ \ \ \ } \end{array}
 \\
 \begin{array}{c}
 \begin{array}{c}
 \begin{array}{ccc}
{{ a}}_{{ p}{ -}{ q}{ +}{ 2}} &  \begin{array}{ccc}
{{ a}}_{{ p}{ -}{ q}{ +}{ 3}} & \cdots  \end{array}
 & {{ a}}_{{ p}{ +}{ 1}{ \ }} \end{array}
 \end{array}
 \\
 \begin{array}{ccc}
\vdots { \  } &  \begin{array}{ccc}
  &  \ \ \ & \vdots  \end{array}
 &  \begin{array}{ccc}
  &  \ \ \  & \vdots  \end{array}
 \end{array}
 \end{array}
 \\
 \begin{array}{ccc}
{{ \ \ \ }{ a}}_{{ p}} &  \begin{array}{ccc}
{{ \ \ \ \ \ \ \ \ }{ a}}_{{ p}{ +}{ 1}} & \cdots  \end{array}
 & { \ \ }{{ a}}_{{ p}{ +}{ q}{ -}{ 1}} \end{array}
 \end{array}
\right)}_{{ q}}\ne { 0}\]
where ${{ a}}_{{ i}}{ =}{ 0}$ for${ \ }{ i}<{ 0}$. Whenever the latter holds, we write
\[{ f}\in {\mathfrak{D}}_{{ p},{ q}}\left({ \zeta}\right).\]
In this case, the Pad\'e approximant
\[{\left[{{ p}}/{{ q}}\right]}_{{ f},{ \zeta}}\left({ z}\right){ =}\frac{{\mathcal A}\left({ f},{ \zeta}\right)\left({ z}\right)}{{\mathcal B}\left({ f},{ \zeta}\right)\left({ z}\right)}\]
is given explicitly by the following Jacobi formula (\cite{Bak2}):

\[{\mathcal A}\left({ f},{ \zeta}\right)\left({ z}\right){ =}{ det}\left( \begin{array}{ccc}
{\left({ z}{ -}{ \zeta}\right)}^{{ q}}{{ S}}_{{ p}{ -}{ q}}\left({ f},{ \zeta}\right)\left({ z}\right) & {\left({ z}{ -}{ \zeta}\right)}^{{ q}{ -}{ 1}}{{ S}}_{{ p}{ -}{ q}{ +}{ 1}}\left({ f},{ \zeta}\right)\left({ z}\right) &  \begin{array}{cc}
{ \dots } & {{ S}}_{{ p}}\left({ f},{ \zeta}\right)\left({ z}\right) \end{array}

 \\
{{ a}}_{{ p}{ -}{ q}{ +}{ 1}} & {{ a}}_{{ p}{ -}{ q}{ +}{ 2}} &  \begin{array}{cc}
{ \dots } &  {{ a}}_{{ p}{ +}{ 1}{ \ }} \end{array}
 \\
 \begin{array}{c}
\vdots  \\
{{ \ }{ a}}_{{ p}} \end{array}
 &  \begin{array}{c}
\vdots  \\
{{ \ \ \ \ \ }{ a}}_{{ p}{ +}{ 1}}{ \ \ } \end{array}
 &  \begin{array}{cc}
 \begin{array}{c}
 \\
\cdots  \end{array}
 &  \begin{array}{c}
\vdots  \\
{ \ }{{ a}}_{{ p}{ +}{ q}} \end{array}
 \end{array}
 \end{array}
\right)\ ,\]
\[{\mathcal B}\left({ f},{ \zeta}\right)\left({ z}\right){ =}{ det}\left( \begin{array}{ccc}
({{ z}-\zeta})^{{ q}} & ({{ z}-\zeta})^{{ q}{ -}{ 1}} &  \begin{array}{cc}
{ \dots \ } & { \ \ }{ 1}{ \ \ \ \ \ } \end{array}
 \\
{{ a}}_{{ p}{ -}{ q}{ +}{ 1}} & {{ a}}_{{ p}{ -}{ q}{ +}{ 2}} &  \begin{array}{cc}
{ \dots } & {{ a}}_{{ p}{ +}{ 1}{ \ }} \end{array}
 \\
 \begin{array}{c}
\vdots  \\
{{ \ \ \ }{ a}}_{{ p}} \end{array}
 &  \begin{array}{c}
\vdots  \\
{{ \ \ \ \ \ \ \ \ }{ a}}_{{ p}{ +}{ 1}}{ \ \ } \end{array}
 &  \begin{array}{cc}
 \begin{array}{c}
 \\
\cdots  \end{array}
 &  \begin{array}{c}
\vdots  \\
{ \ }{{ a}}_{{ p}{ +}{ q}} \end{array}
 \end{array}
 \end{array}
\right)\]
with
\begin{center}
\[{{ S}}_{{ k}}\left({ f},{ \zeta}\right)\left({ z}\right){ =}\left\{ \begin{array}{c}
\sum^{{ k}}_{{ v}{ =}{ 0}}{{{ a}}_{{ v}}{\left({ z}{ -}{ \zeta}\right)}^{{ v}}}{ ,\ \ }{ if}{ \ }{ k}\ge { 0} \\
{ 0}{ ,\ \ \ }{ if}{ \ \ }{ k}<{ 0}. \end{array}
\right.\]
\end{center}
\begin{def}\label{def2.1}
If${ \ }{\mathcal A}\left({ f},{ \zeta}\right)\left({ z}\right)$ and${ \ }{\mathcal B}\left({ f},{ \zeta}\right)\left({ z}\right)$ are given by the previous Jacobi formula and
${ f}\in {\mathfrak{D}}_{{ p},{ q}}\left({ \zeta}\right)$, then
they do not have a common zero on $\mathbb{C}$.
\end{def}
If ${ K}\subset\mathbb{C}$ is compact and ${ f}\in {\mathfrak{D}}_{{ p},{ q}}\left({ \zeta}\right)$,
then there exists $\delta=\delta({ f},\zeta,{ K})>0$ such that

\begin{center}
 ${\left|{\mathcal A}\left({ f},{ \zeta}\right)\left({ z}\right)\right|}^{{ 2}}{ +}{\left|{\mathcal B}\left({ f},{ \zeta}\right)\left({ z}\right)\right|}^{{ 2}}{ >}\delta$ for all ${ z}\in { K}$.
 \end{center}

\noindent We shall also make use of the following standard proposition, regarding the Pad\'e approximants of rational functions (\cite{Bak2}, Theorem 1.4.4).

\begin{prop}\label{prop2.2}
 Let$\ \ \Phi\left({ z}\right){ =}{{\mathcal A}\left({ z}\right)}/{{\mathcal B}\left({ z}\right)}$ be a rational function, where ${\mathcal A}$ and ${\mathcal B}$ are polynomials without any common zero in$\ {\mathbb C}$. Let ${ deg}{\mathcal A}{ =}{ \kappa}$ and$\ \ { deg}{\mathcal B}{ =}{ \lambda}$. Then, for every ${ \zeta}\in {\mathbb C}$ such that ${\mathcal B}\left({ \zeta}\right)\ne { 0}$, the following hold
\begin{center}
$\Phi\in {\mathfrak{D}}_{{ \kappa},{ \lambda}}\left({ \zeta}\right)$,\\
$\Phi\in {\mathfrak{D}}_{{ p},{ \lambda}}\left({ \zeta}\right)$ for all ${ p}>{ \kappa}$,\\
$\Phi \in {\mathfrak{D}}_{{ \kappa},{ q}}\left({ \zeta}\right)$ for all ${ q}>{ \lambda}$.
\end{center}
Furthermore, if ${ p}>{ \kappa}$ and$\ \ { q}>{ \lambda}$, then
\[\ \ \Phi\notin {\mathfrak{D}}_{{ p},{ q}}\left({ \zeta}\right).\]

\noindent
In all these cases $\Phi$ coincides with its corresponding Pad\'e approximant.

\end{prop}
\section{UNIVERSAL PAD\'E APPROXIMANTS IN ${\mathcal H}\left({ \Omega}\right)$}
\noindent

We start by proving the following theorem.

\begin{thm}\label{thm3.1}
Let$\ \ \mathfrak{F}$ be a subset of ${\mathbb N}\times {\mathbb N}$ containing a sequence$\ \ \left({\widetilde{{ p}}}_{{ n}},{\widetilde{{ q}}}_{{ n}}\right)\in \mathfrak{F}$, ${ n}{ =}{ 1},{ 2}{ ,\dots }$, with ${\widetilde{{ p}}}_{{ n}}\to { +}\infty $. Let ${ \Omega}\subset {\mathbb C}$ be a simply connected open set and ${ L}$ and ${ L}'$ two compact subsets of ${ \Omega}$ (with connected complements). Let ${ K}\subset {\mathbb C}\backslash { \Omega}$ be another compact set with connected complement. Then there exists a holomorphic function$\ { f}\in {\mathcal H}\left({ \Omega}\right)$, such that for every polynomial ${ h}$ there exists a sequence $\left({{ p}}_{{ n}},{{ q}}_{{ n}}\right)\in \mathfrak{F}$, ${ n}{ =}{ 1},{ 2}{ ,\dots }$, such that the following hold:
\begin{enumerate}
\item  ${ f}\in {\mathfrak{D}}_{{{ p}}_{{ n}},{{ q}}_{{ n}}}\left({ \zeta}\right)$
for every ${ n}\in {\mathbb N}$ and ${ \zeta}\in { L}$,

\item  ${{ sup}}_{{ \zeta}\in { L}{\rm ,\ }{ z}\in { K}}\left|{\left[{{{ p}}_{{ n}}}/{{{ q}}_{{ n}}}\right]}_{{ f},{ \zeta}}\left({ z}\right){ -}{ h}\left({ z}\right)\right|\to { 0}$, as ${ n}\to { +}\infty $ and

\item  ${{ sup}}_{{ \zeta}\in { L}{\rm ,\ }{ z}\in { L}'}\left|{\left[{{{ p}}_{{ n}}}/{{{ q}}_{{ n}}}\right]}_{{ f},{ \zeta}}\left({ z}\right){ -}{ f}\left({ z}\right)\right|\to { 0}$, as$\ \ { n}\to { +}\infty $.
\end{enumerate}
The set of all such functions ${ f}\in {\mathcal H}\left({ \Omega}\right)$ is dense and ${{ G}}_{{ \delta}}$ in$\ \ {\mathcal H}\left({ \Omega}\right)$, endowed with the topology of uniform convergence on compacta.
\end{thm}
\begin{Proof} Let ${{ f}}_{{ j}}$, ${ j}{ =}{ 1},{ 2}{ ,\dots }$ be an enumeration of all polynomials with coefficients in ${\mathbb Q}{ +}{ i}{\mathbb Q}$. The set of functions ${ f}\in {\mathcal H}\left({ \Omega}\right)$ satisfying Theorem 3.1 can be described as follows:
\[{\mathcal Y}\left({ \Omega}{ ,\ }{{ L}}^{{ '}}{ ,\ }{ L},{ K}\right){ :=}\bigcap^{\infty }_{{ j},{ s}{ =}{ 1}}{\bigcup_{\left({ p},{ q}\right)\in \mathfrak{F}}{\left[{\mathcal G}\left({ \Omega}{ ,\ }{{ L}}^{{ '}},{ L},{ p},{ q},{ s}\right)\bigcap {\mathcal F}\left({ \Omega},{ K},{ L},{ p},{ q},{ j},{ s}\right)\right]}}\]
where
\[{\mathcal G}\left({ \Omega}{ ,\ }{{ L}}^{{ '}},{ L},{ p},{ q},{ s}\right){ =}\left\{{ f}\in {\mathcal H}\left({ \Omega}\right){ :}{ f}\in {\mathfrak{D}}_{{ p},{ q}}\left({ \zeta}\right){ \ }{ for}{ \ }{ all}{ \ }{ \zeta}\in { L}\right.{ \  }{ and}\]
\[\left.{{ sup}}_{{ \zeta}\in { L}{ ,\ }{ z}\in {{ L}}^{{ '}}}\left|{\left[{{ p}}/{{ q}}\right]}_{{ f},{ \zeta}}\left({ z}\right){ -}{ f}\left({ z}\right)\right|{ <}\frac{{ 1}}{{ s}}\right\}\]
and
\[{\mathcal F}\left({ \Omega},{ K},{ L},{ p},{ q},{ j},{ s}\right){ =}\left\{{ f}\in {\mathcal H}\left({ \Omega}\right){ :}{ f}\in {\mathfrak{D}}_{{ p},{ q}}\left({ \zeta}\right){ \ }{ for}{ \ }{ all}{ \ }{ \zeta}\in { L}\right.{\ and}\]
\[\left.{{ sup}}_{{ \zeta}\in { L}{ ,\ }{ z}\in { K}}\left|{\left[{{ p}}/{{ q}}\right]}_{{ f},{ \zeta}}\left({ z}\right){ -}{{ f}}_{{ j}}\left({ z}\right)\right|{ <}\frac{{ 1}}{{ s}}\right\}\]

First we show that the sets ${\mathcal G}\left({ \Omega}{ ,\ }{{ L}}^{{ '}},{ L},{ p},{ q},{ s}\right)$ and ${\mathcal F}\left({ \Omega},{ K},{ L},{ p},{ q},{ j},{ s}\right)$ are open in ${\mathcal H}\left({ \Omega}\right)$ which imply that ${\mathcal Y}\left({ \Omega}{ ,\ }{{ L}}^{{ '}}{ ,\ }{ L},{ K}\right)$ is ${{ G}}_{{ \delta}}$ in$\ \ {\mathcal H}\left({ \Omega}\right)$.

Let$\ \ { f}\in {\mathcal G}\left({ \Omega}{ ,\ }{{ L}}^{{ '}},{ L},{ p},{ q},{ s}\right)$. Let also $\widetilde{{ L}}\subset { \Omega}$ be a compact set such that$\ \ { L}\bigcup {{ L}}^{{ '}}\subset {\left(\widetilde{{ L}}\right)}^{{ o}}$. There exists an$\ \ { r}>0$, such that $\left\{{ z}\in {\mathbb C}{\rm :}\left|{ z}{ -}{ \zeta}\right|\le { r}\right\}\subset {\left(\widetilde{{ L}}\right)}^{{ o}}$ for every$\ \ { \zeta}\in { L}$. Let us consider ${ a}>0$ small and a function ${ g}\in {\mathcal H}\left({ \Omega}\right)$ satisfying$\ \ {{ sup}}_{{ z}\in \widetilde{{ L}}}\left|{ f}\left({ z}\right){ -}{ g}\left({ z}\right)\right|{ <}{ a}$. We will show that if ${ a}$ is small enough, then ${ g}\in {\mathcal G}\left({ \Omega}{ ,\ }{{ L}}^{{ '}},{ L},{ p},{ q},{ s}\right)$.

By the Cauchy estimates, we can uniformly control any finite set of Taylor coefficients of${ \ \ }{ g}$ with center$\ \ { \zeta}\in { L}$, to be close to the corresponding coefficients of$\ \ { f}$.

The Hankel determinant defining${{ \ \ }\mathfrak{D}}_{{ p},{ q}}\left({ \zeta}\right)$ for ${ f}$ depends continuously on$\  { \zeta}\in { L}$; thus, its absolute value is greater than a$\ \ { \delta}>0$. If we choose ${ a}>0$ small enough, then the absolute values of the Hankel determinants defining ${{ \ \ }\mathfrak{D}}_{{ p},{ q}}\left({ \zeta}\right)$ for ${ g}$ will be greater than ${{ \delta}}/{{ 2}}$ for all$\ \ { \zeta}\in { L}$. Therefore ${ g}\in {\mathfrak{D}}_{{ p},{ q}}\left({ \zeta}\right)$, for all$\ \ { \zeta}\in { L}$. We recall the Jacobi formulas
\begin{align}\label{3.1}
[p/q]_{f,\zeta}( z) =\frac{{\mathcal A}( f,\zeta)( z)}{\mathcal{B}( f,\zeta)( z)}
\end{align}
and
\begin{align}\label{3.2}
{\left[{{ p}}/{{ q}}\right]}_{{ g},{ \zeta}}\left({ z}\right){ =}\frac{{\mathcal A}\left({ g},{ \zeta}\right)\left({ z}\right)}{{\mathcal B}\left({ g},{ \zeta}\right)\left({ z}\right)}.
\end{align}
The coefficients of the polynomials ${\mathcal A}\left({ f},{ \zeta}\right)\left({ z}\right)$ and ${\mathcal B}\left({ f},{ \zeta}\right)\left({ z}\right)$ vary continuously with${ \ \ }{ \zeta}\in { L}$. Since ${\mathcal A}\left({ f},{ \zeta}\right)\left({ z}\right)$ and ${\mathcal B}\left({ f},{ \zeta}\right)\left({ z}\right)$ do not have a common zero in $\mathbb{C}$, there exists ${ \delta}{ '}>0$ such that
\begin{align}\label{3.3}
{\left|{\mathcal A}\left({ f},{ \zeta}\right)\left({ z}\right)\right|}^{{ 2}}{ +}{\left|{\mathcal B}\left({ f},{ \zeta}\right)\left({ z}\right)\right|}^{{ 2}}{ >}{{ \delta}}^{{ '}}
\end{align}
for all ${ \zeta}\in { L}$ and ${ z}\in { L}{ '}$. If ${ a}>0$ is small enough, then ${\left|{\mathcal A}\left({ g},{ \zeta}\right)\left({ z}\right)\right|}^{{ 2}}{ +}{\left|{\mathcal B}\left({ g},{ \zeta}\right)\left({ z}\right)\right|}^{{ 2}}{ >}{{{ \delta}}^{{ '}}}/{{ 2}}$ for all ${ \zeta}\in { L}$ and ${ z}\in { L}{ '}$. 
We set
\begin{align}\label{3.4}
{ \gamma}{ =}\frac{{ 1}}{{ s}}{ -}{{ sup}}_{{ \zeta}\in { L}{ ,\ }{ z}\in {{ L}}^{{ '}}}\left|{\left[{{ p}}/{{ q}}\right]}_{{ f},{ \zeta}}\left({ z}\right){ -}{ f}\left({ z}\right)\right|
\end{align}
and we notice that ${ \gamma}>0$. If ${ a}>0$ is sufficiently small, then
\begin{align}\label{3.5}
{{ sup}}_{{ \ }{ z}\in {{ L}}^{{ '}}}\left|{ f}\left({ z}\right){ -}{ g}\left({ z}\right)\right|{ <}\frac{{ \gamma}}{{ 2}}.
\end{align}
Since ${\left[{{ p}}/{{ q}}\right]}_{{ f},{ \zeta}}\left({ z}\right)$ takes only finite  values for ${ \zeta}\in { L}$ and$\ \ { z}\in {{ L}}^{{ '}}$, it is uniformly bounded by a constant, say$\ \ { M}<+\infty $. Relation ~\eqref{3.1} and ~\eqref{3.3} imply
\begin{align*}
\frac{{{ \delta}}^{{ '}}}{{\left|{\mathcal B}\left({ f},{ \zeta}\right)\left({ z}\right)\right|}^{{ 2}}}{ <}{{ M}}^{{ 2}}{ +}{ 1}.
\end{align*}
for all ${ \zeta}\in { L}$ and ${ z}\in { L}{ '}$. Hence, there exists a ${{ \delta''}}{ >}0$ so that
\begin{align}\label{3.6}
{{ \delta''}}{ <}\left|{\mathcal B}\left({ f},{ \zeta}\right)\left({ z}\right)\right|
\end{align}
for all ${ \zeta}\in { L}$ and ${ z}\in { L}{ '}$. By continuous
dependence, if ${ a}>{ 0}$ is small, we have
\begin{align}\label{3.7}
\frac{\delta''}{2}<|\mathcal{B}(g,\zeta)(z)|
\end{align}
for all ${ \zeta}\in { L}$ and ${ z}\in { L}{ '}$. Relation ~\eqref{3.1}, ~\eqref{3.2}, ~\eqref{3.6} and ~\eqref{3.7} imply
\begin{align}\label{3.8}
{{ sup}}_{{ \zeta}\in { L}{ ,\ }{ z}\in {{ L}}^{{ '}}}\left|{\left[{{ p}}/{{ q}}\right]}_{{ g},{ \zeta}}\left({ z}\right){ -}{\left[{{ p}}/{{ q}}\right]}_{{ f},{ \zeta}}\left({ z}\right)\right|{ <}\frac{{ \gamma}}{{ 2}}
\end{align}
provided that ${ a}>0$ is small enough. Relations ~\eqref{3.4}, ~\eqref{3.5}, ~\eqref{3.8} combined with the triangle inequality imply
$${{ sup}}_{{ \zeta}\in { L}{ ,\ }{ z}\in {{ L}}^{{ '}}}\left|{\left[{{ p}}/{{ q}}\right]}_{{ g},{ \zeta}}\left({ z}\right){ -}{ g}\left({ z}\right)\right|{ <}\frac{{ 1}}{{ s}},$$
which shows that the set ${\mathcal G}\left({ \Omega}{ ,\ }{{ L}}^{{ '}},{ L},{ p},{ q},{ s}\right)$ is open in ${\mathcal H}\left({ \Omega}\right)$.

Next we prove that the set ${\mathcal F}\left({ \Omega},{ K},{ L},{ p},{ q},{ j},{ s}\right)$ is open in$\ \ {\mathcal H}\left({ \Omega}\right)$. Let ${ f}\in {\mathcal F}\left({ \Omega},{ K},{ L},{ p},{ q},{ j},{ s}\right)$ and $\widetilde{{ L}}\subset { \Omega}$ be a compact set such that$\ \ { L}\subset {\left(\widetilde{{ L}}\right)}^{{ o}}$. We consider $ { a}>0$ small and a function ${ g}\in {\mathcal H}\left({ \Omega}\right)$ satisfying$\ \ {{ sup}}_{{ \ }{ z}\in \widetilde{{ L}}}\left|{ f}\left({ z}\right){ -}{ g}\left({ z}\right)\right|{ <}{ a}$. We will see that if ${ a}>0$ is small enough, then ${ g}\in {\mathcal F}\left({ \Omega},{ K},{ L},{ p},{ q},{ j},{ s}\right)$.

By the Cauchy estimates we can uniformly control any finite set of Taylor coefficients of ${ g}$ with center ${ \zeta}\in { L}$, to be close to the corresponding coefficients of ${ f}$. The proofs that ${ g}\in {\mathfrak{D}}_{{ p},{ q}}\left({ \zeta}\right)$ for all ${ \zeta}\in { L}$ are similar to the previous ones and are omitted. It remains to control the quantity
\[{{ sup}}_{{ \zeta}\in { L}{ ,\ }{ z}\in { K}}\left|{\left[{{ p}}/{{ q}}\right]}_{{ g},{ \zeta}}\left({ z}\right){ -}{\left[{{ p}}/{{ q}}\right]}_{{ f},{ \zeta}}\left({ z}\right)\right|.\]

Again by the Cauchy estimates and the Jacobi formula the coefficients of the polynomials$\ \ {\mathcal A}\left({ f},{ \zeta}\right)\left({ z}\right)$, ${\mathcal B}\left({ f},{ \zeta}\right)(z)$, ${\mathcal A}\left({ g},{ \zeta}\right)\left({ z}\right)$ and ${\mathcal B}\left({ g},{ \zeta}\right)(z)$ vary continuously. Since
${\mathcal A}\left({ f},{ \zeta}\right)\left({ z}\right)$ and ${\mathcal B}\left({ f},{ \zeta}\right)\left({ z}\right)$ do not have a common zero in $\mathbb{C}$
and ${{ f}}_{{ j}}\left({ z}\right)$ does not have any poles in$\ \ {\mathbb C}$, the inequality
\[{{ sup}}_{{ \zeta}\in { L}{ ,\ }{ z}\in { K}}\left|{\left[{{ p}}/{{ q}}\right]}_{{ f},{ \zeta}}\left({ z}\right){ -}{{ f}}_{{ j}}\left({ z}\right)\right|{ <}\frac{{ 1}}{{ s}}\]
implies that the approximant ${\left[{{ p}}/{{ q}}\right]}_{{ f},{ \zeta}}\left({ z}\right){ =}{{\mathcal A}\left({ f},{ \zeta}\right)\left({ z}\right)}/{{\mathcal B}\left({ f},{ \zeta}\right)}$ is such that $ \inf_{{ \zeta}\in { L}{ ,\ }{ z}\in { K}}\left|{\mathcal B}\left({ f},{ \zeta}\right)\right|{ >}0$. It follows that
\[{{ sup}}_{{ \zeta}\in { L}{ ,\ }{ z}\in { K}}\left|{\left[{{ p}}/{{ q}}\right]}_{{ g},{ \zeta}}\left({ z}\right){ -}{\left[{{ p}}/{{ q}}\right]}_{{ f},{ \zeta}}\left({ z}\right)\right|\]
is small. The triangle inequality implies
\[{{ sup}}_{{ \zeta}\in { L}{ ,\ }{ z}\in { K}}\left|{\left[{{ p}}/{{ q}}\right]}_{{ g},{ \zeta}}\left({ z}\right){ -}{{ f}}_{{ j}}\left({ z}\right)\right|{ <}\frac{{ 1}}{{ s}}, \]
provided that ${ a}>0$ is small enough. Thus, ${\mathcal F}\left({ \Omega},{ K},{ L},{ p},{ q},{ j},{ s}\right)$ is open in$\ \ {\mathcal H}\left({ \Omega}\right)$.

It remains to prove that the set
\[\bigcup_{\left({ p},{ q}\right)\in \mathfrak{F}}{{ \ }\left[{\mathcal G}\left({ \Omega}{ ,\ }{{ L}}^{{ '}},{ L},{ p},{ q},{ s}\right)\bigcap {\mathcal F}\left({ \Omega},{ K},{ L},{ p},{ q},{ j},{ s}\right)\right]}\]
is dense in ${\mathcal H}\left({ \Omega}\right)$. Let$\ \ {{ L''}}\subset { \Omega}$ be a compact set, ${ \Phi}\in {\mathcal H}\left({ \Omega}\right)$ and $\varepsilon>0$. We are looking for a ${ g}\in {\mathcal H}\left({ \Omega}\right)$ and a $\left({ p},{ q}\right)\in \mathfrak{F}$ such that
\noindent ${ g}\in {\mathfrak{D}}_{{ p},{ q}}\left({ \zeta}\right)$ for all ${ \zeta}\in { L}$,
\[{{ sup}}_{{ \ }{ z}\in {{ L''}}}\left|{ g}\left({ z}\right){ -}{ \Phi}\left({ z}\right)\right| <\varepsilon, \]
${{ sup}}_{{ \zeta}\in { L}{ ,\ }{ z}\in {{ L'}}}\left|{\left[{{ p}}/{{ q}}\right]}_{{ g},{ \zeta}}\left({ z}\right){ -}{ g}\left({ z}\right)\right|{ <}\frac{{ 1}}{{ s}}$ and
\[{{ sup}}_{{ \zeta}\in { L}{ ,\ }{ z}\in { K}}\left|{\left[{{ p}}/{{ q}}\right]}_{{ g},{ \zeta}}\left({ z}\right){ -}{{ f}}_{{ j}}\left({ z}\right)\right|{ <}\frac{{ 1}}{{ s}}. \]
Without loss of generality, we assume that ${ L}\bigcup {{ L}}^{{ '}}\subset {\left({{ L''}}\right)}^{{ o}}$ and that the compact set ${{ L''}}$ has connected complement.

We consider the function
\[{\mathcal W}\left({ z}\right){ =}\left\{ \begin{array}{c}
{{ f}}_{{ j}}\left({ z}\right){ ,\ }{ if}{ \ }{ z}\in { K} \\
{ \Phi}\left({ z}\right){ ,\ }{ if}{ \ }{ z}\in {{ L''}} \end{array}
\right.. \]
The function ${\mathcal W}$ is well defined because$\ \ { K}\bigcap {{ L''}}{ =}\emptyset $. Since ${ K}\bigcup {{ L''}}$ is a compact set with connected complement, the function ${\mathcal W}$ can be approximated by a polynomial ${ P}\left({ z}\right)$ uniformly on$\ \ { K}\bigcup {{ L''}}$.
Our assumption on $\mathfrak{F}$ allows us to find $\left({ p},{ q}\right)\in \mathfrak{F}$ such that$\ \ { p}>deg{P}$.

We consider the function $\mathfrak{A}\left({ z}\right){ =}{ P}\left({ z}\right){ +}{ d}{{ z}}^{{ p}}$ where ${ d}\ne { 0}$ is small in absolute value. Then, $\ \mathfrak{A}\left({ z}\right)$ and$\ \ {\mathcal W}\left({ z}\right)$ are uniformly close on$\ \ { K}\bigcup {{ L''}}$. It follows that

\noindent $\mathfrak{A}\left({ z}\right)\in {\mathfrak{D}}_{{ p},{ q}}\left({ \zeta}\right)$ and ${\left[{{ p}}/{{ q}}\right]}_{\mathfrak{A},{ \zeta}}\left({ z}\right)\equiv \mathfrak{A}\left({ z}\right)$
 for all ${ z}\in {\mathbb C}$; in particular, this holds for all ${ z}\in { L}$.
Thus, the quantity
\[{{ sup}}_{{ \zeta}\in { L}{ ,\ }{ z}\in { K}}\left|{\left[{{ p}}/{{ q}}\right]}_{\mathfrak{A},{ \zeta}}\left({ z}\right){ -}{{ f}}_{{ j}}\left({ z}\right)\right|\]
is small. Moreover
\[{{ sup}}_{{ \zeta}\in { L}{ ,\ }{ z}\in {{ L}}^{{ '}}}\left|{\left[{{ p}}/{{ q}}\right]}_{\mathfrak{A},{ \zeta}}\left({ z}\right){ -}\mathfrak{A}\left({ z}\right)\right|{ =}{ 0}<\frac{{ 1}}{{ s}}.\]
We also note that the quantity
\[{{ sup}}_{{ z}\in {{ L''}}}\left|{\left[{{ p}}/{{ q}}\right]}_{\mathfrak{A},{ \zeta}}\left({ z}\right){ -}{ \Phi}\left({ z}\right)\right|\]
is small, provided that ${ 0}<\left|{ d}\right|$ is small enough. We set
\[{ g}\left({ z}\right)\equiv \mathfrak{A}\left({ z}\right){ =}{ P}\left({ z}\right){ +}{ d}{{ z}}^{{ p}}\in {\mathcal H}\left({ \Omega}\right)\]
and the proof is completed. Baire's category theorem yields the result.$\qquad\blacksquare$
\end{Proof}
We now recall the following lemma (\cite{Nest2}).
\begin{lem}\label{lem3.2}
Let$\ \ { \Omega}$ be a domain in$\ {\mathbb C}$. Then there exists a sequence${{ \ \ }{ K}}_{{ m}}$, ${ m}{ =}{ 1},{ 2}{ ,\dots }$, of compact subsets of ${\mathbb C}\backslash { \Omega}$ with connected complements, such that for every compact set ${ K}\subset {\mathbb C}\backslash { \Omega}$ with connected complement,
there exists an ${ m}\in \left\{{ 1},{ 2}{ ,\dots }\right\}$ so that ${ K}\subset {{ K}}_{{ m}}$. ¦
\end{lem}
If we set ${ L}{ =}\left\{{ \zeta}\right\}$ to be a singleton and ${ K}{ =}{{ K}}_{{ m}}$ be given by Lemma 3.2 and ${{ L}}^{{ '}}{ =}{{ L}}^{{ '}}_{{ n}}$ to vary in an exhausting sequence of compact sets of$\ \ { \Omega}$, then applying Baire's theorem we obtain the following theorem.

\begin{thm}\label{thm3.3}
Let$\ \ \mathfrak{F}$ be a subset of ${\mathbb N}\times {\mathbb N}$ containing a sequence$\ \ \left({\widetilde{{ p}}}_{{ n}},{\widetilde{{ q}}}_{{ n}}\right)\in \mathfrak{F}$,${ \ }{ n}{ =}{ 1},{ 2}{ ,\dots }$, with ${\widetilde{{ p}}}_{{ n}}\to { +}\infty $. Let ${ \Omega}\subset {\mathbb C}$ be a simply connected domain and ${ \zeta}\in { \Omega}$ be fixed. Then there exists a holomorphic function$\ \ { f}\in {\mathcal H}\left({ \Omega}\right)$, such that, for every polynomial ${ h}$ and every compact set ${ K}\subset {\mathbb C}\backslash { \Omega}$ with connected complement, there exists a sequence $\ \ \left({{ p}}_{{ n}},{{ q}}_{{ n}}\right)$ $\in \mathfrak{F}$,${ \ \ }{ n}{ =}{ 1},{ 2}{ ,\dots }$, such that the following hold:
\begin{enumerate}
\item  ${ f}\in {\mathfrak{D}}_{{{ p}}_{{ n}},{{ q}}_{{ n}}}\left({ \zeta}\right)$ for all ${ n}\in {\mathbb N}$,

\item  ${{ sup}}_{{\rm \ }{ z}\in { K}}\left|{\left[{{{ p}}_{{ n}}}/{{{ q}}_{{ n}}}\right]}_{{ f},{ \zeta}}\left({ z}\right){ -}{ h}\left({ z}\right)\right|\to { 0}$, as ${ n}\to { +}\infty $ and

\item  ${{ sup}}_{{ z}\in { L}'}\left|{\left[{{{ p}}_{{ n}}}/{{{ q}}_{{ n}}}\right]}_{{ f},{ \zeta}}\left({ z}\right){ -}{ f}\left({ z}\right)\right|\to { 0}$, as$\ \ { n}\to { +}\infty $.
\end{enumerate}
The set of all such functions ${ f}\in {\mathcal H}\left({ \Omega}\right)$ is dense and ${{ G}}_{{ \delta}}$ in$\ \ {\mathcal H}\left({ \Omega}\right)$.
\end{thm}
Setting ${ L}{ =}{{ L}}^{{ '}}{ =}{{ L}}_{{ n}}$ to be an exhausting sequence of compact subsets of ${ \Omega}$ and ${ K}{ =}{{ K}}_{{ m}}$ to be given by Lemma 3.2 and applying Baire's category theorem once more we obtain the following.

\begin{thm}\label{thm3.4}
Let$\ \ \mathfrak{F}$ be a subset of ${\mathbb N}\times {\mathbb N}$ containing a sequence$\ \ \left({\widetilde{{ p}}}_{{ n}},{\widetilde{{ q}}}_{{ n}}\right)\in \mathfrak{F}$,${ \ }{ n}{ =}{ 1},{ 2}{ ,\dots }$, with ${\widetilde{{ p}}}_{{ n}}\to { +}\infty $. Let also ${ \Omega}\subset {\mathbb C}$ be a simply connected domain. Then there exists a holomorphic function$\ \ { f}\in {\mathcal H}\left({ \Omega}\right)$ satisfying the following.
For every compact set ${ K}\subset {\mathbb C}\backslash { \Omega}$ with connected complement and every polynomial$\ \ { h}$, there exists a sequence$\ \ \left({{ p}}_{{ n}},{{ q}}_{{ n}}\right)\in \mathfrak{F}$,${ \ \ }{ n}{ =}{ 1},{ 2}{ ,\dots }$, such that
for every compact set$\ \ { L}\subset { \Omega}$ we have

\begin{enumerate}
\item  ${ f}\in {\mathfrak{D}}_{{{ p}}_{{ n}},{{ q}}_{{ n}}}\left({ \zeta}\right)$ for all ${ n}\in {\mathbb N}$ and ${ \zeta}\in { L}$,

\item  ${{ sup}}_{{\rm \ }{ \zeta}\in { L}{ ,\ }{ z}\in { K}}\left|{\left[{{{ p}}_{{ n}}}/{{{ q}}_{{ n}}}\right]}_{{ f},{ \zeta}}\left({ z}\right){ -}{ h}\left({ z}\right)\right|\to { 0}$, as ${ n}\to { +}\infty $ and

\item  ${{ sup}}_{{ \zeta}\in { L}{ ,\ }{ z}\in { L}}\left|{\left[{{{ p}}_{{ n}}}/{{{ q}}_{{ n}}}\right]}_{{ f},{ \zeta}}\left({ z}\right){ -}{ f}\left({ z}\right)\right|\to { 0}$, as$\ \ { n}\to { +}\infty $.
\end{enumerate}
The set of all such functions ${ f}\in {\mathcal H}\left({ \Omega}\right)$ is dense and ${{ G}}_{{ \delta}}$ in$\ \ {\mathcal H}\left({ \Omega}\right)$. ¦
\end{thm}
\begin{rem}\label{rem3.5}
If, in Theorems 3.3 and 3.4, we consider only compact sets ${ K}\subset {\mathbb C}\backslash \overline{{ \Omega}}$ having connected complements, then the corresponding classes are also dense and ${{ G}}_{{ \delta}}$ in$\ \ {\mathcal H}\left({ \Omega}\right)$. These classes generalize the classes of universal Taylor series in the sense of Luh and Chui and Parnes, but there is no guarantee that the universal functions extend
continuously on the closure of Omega. This will be obtained in section 4
below.
\end{rem}
\begin{rem}\label{rem3.6}
Since the compact set ${ K}$ has connected complement, by Mergelyan's theorem, in the previous results
the function ${ h}:{ K}\to\mathbb{C}$ can be more generally assumed to be in $A({ K})$, that is,
continuous in ${ K}$ and holomorphic in the interior ${ K}^0$.
\end{rem}

\section{UNIVERSALITY IN $A\left({\rm \Omega}\right)$ }

Let ${ \Omega}\subset {\mathbb C}$ be a domain, such that ${\overline{{ \Omega}}}^{{ o}}{ =}{ \Omega}$ and $\left\{\infty \right\}\bigcup \left({\mathbb C}\backslash \overline{{ \Omega}}\right)$ is connected. In what follows the domain ${ \Omega}$ always satisfies these two assumptions. The space ${ A}\left({ \Omega}\right)\equiv {\mathcal H}\left({ \Omega}\right)\bigcap { C}\left(\overline{{ \Omega}}\right)$ is a Frechet space with the seminorms ${{ sup}}_{{ z}\in \overline{{ \Omega}}{ ,\ }\left|{ z}\right|\le { n}}\left|{ f}\left({ z}\right)\right|$, ${ f}\in { A}\left({ \Omega}\right)$, ${ n}\in\{{ 1},{ 2}{ ,\dots }\}$. We have the following result.
\begin{thm}\label{thm4.1}
Let$\ \ \mathfrak{F}$ be a subset of ${\mathbb N}\times {\mathbb N}$ containing a sequence$\ \ \left({\widetilde{{ p}}}_{{ n}},{\widetilde{{ q}}}_{{ n}}\right)\in \mathfrak{F}$, ${ n}{ =}{ 1},{ 2}{ ,\dots }$, with ${\widetilde{{ p}}}_{{ n}}\to { +}\infty $. Let ${ \Omega}\subset {\mathbb C}$ be as above and ${ L}\subset { \Omega}$, ${ L}'\subset \overline{{ \Omega}}$ two compact subsets with connected complements. Let also ${K}$ be another compact set
 with connected complement such that ${ K}\bigcap \overline{{ \Omega}}{ =}\emptyset $. Then there exists a function ${ f}\in { A}\left({ \Omega}\right)$, such that for every polynomial ${ P}$ there exists a sequence $\left({{ p}}_{{ n}},{{ q}}_{{ n}}\right)\in \mathfrak{F}$, ${ n}{ =}{ 1},{ 2}{ ,\dots }$, such that the following hold:
\begin{enumerate}
\item  ${ f}\in {\mathfrak{D}}_{{{ p}}_{{ n}},{{ q}}_{{ n}}}\left({ \zeta}\right)$
for every ${ n}\in {\mathbb N}$ and ${ \zeta}\in { L}$,

\item  ${{ sup}}_{{ \zeta}\in { L}{\rm ,\ }{ z}\in { K}}\left|{\left[{{{ p}}_{{ n}}}/{{{ q}}_{{ n}}}\right]}_{{ f},{ \zeta}}\left({ z}\right){ -}{ P}\left({ z}\right)\right|\to { 0}$, as ${ n}\to { +}\infty $ and

\item  ${{ sup}}_{{ \zeta}\in { L}{\rm ,\ }{ z}\in { L}'}\left|{\left[{{{ p}}_{{ n}}}/{{{ q}}_{{ n}}}\right]}_{{ f},{ \zeta}}\left({ z}\right){ -}{ f}\left({ z}\right)\right|\to { 0}$, as$\ \ { n}\to { +}\infty $.
\end{enumerate}
The set of all such functions ${ f}\in { A}\left({ \Omega}\right)$ is dense and ${{ G}}_{{ \delta}}$ in$\ \ { A}\left({ \Omega}\right)$.
\end{thm}
\begin{Proof}
The proof is similar to the one of Theorem 3.1 with the obvious modifications.
\end{Proof}
Next, we recall the following lemma (\cite{Luh2}).
\begin{lem}\label{lem4.2}
There exists a sequence${{ \ \ }{ K}}_{{ m}}\subset {\mathbb C}\backslash \overline{{ \Omega}}$,$\ \ \ { m}{ =}{ 1},{ 2}{ ,\dots }$, of compact subsets with connected complements, such that for every compact set ${ K}\subset {\mathbb C}\backslash \overline{{ \Omega}}$ with connected complement,
there exists a ${ m}\in \left\{{ 1},{ 2}{ ,\dots }\right\}$ so that ${ K}\subset {{ K}}_{{ m}}$.
\end{lem}
Varying ${ K}$ to$\ \ {{ K}}_{{ m}}$,$\ \ { m}{ =}{ 1},{ 2}{ ,\dots }$ according to Lemma 4.2, ${ L}'\equiv {{ L}}^{{ '}}_{{ m}}{ =}\left\{{ z}{ :\ }\left|{ z}\right|\le { m}\right\}\bigcap \overline{{ \Omega}}$ and setting$\ \ { L}{ =}\left\{{{ \zeta}}_{{ 0}}\right\}$, ${{ \zeta}}_{{ 0}}\in { \Omega}$, an application of Baire's category theorem yields the following.
\begin{thm}\label{thm4.3}
Let$\ \ \mathfrak{F}$ be a subset of ${\mathbb N}\times {\mathbb N}$ containing a sequence$\ \ \left({\widetilde{{ p}}}_{{ n}},{\widetilde{{ q}}}_{{ n}}\right)\in \mathfrak{F}$,$\ \ { n}{ =}{ 1},{ 2}{ ,\dots }$ with ${\widetilde{{ p}}}_{{ n}}\to { +}\infty $. Let ${ \Omega}$ be a domain in$\ \ {\mathbb C}$, such that $\left\{\infty \right\}\bigcup \left({\mathbb C}\backslash \overline{{ \Omega}}\right)$ is connected and$\ \ {\overline{{ \Omega}}}^{{ o}}{ =}{ \Omega}$. Let ${{ \zeta}}_{{ 0}}$ be fixed. Then there exists a holomorphic function$\ \ { f}\in { A}\left({ \Omega}\right)$, such that, for every polynomial $P({z})$ and any compact set ${ K}\subset {\mathbb C}\backslash \overline{{ \Omega}}$ with connected complement, there exists a sequence$\ \ \left({{ p}}_{{ n}},{{ q}}_{{ n}}\right)\in \mathfrak{F}$,$\ { n}{ =}{ 1},{ 2}{ ,\dots }$, so that
\begin{enumerate}
\item  ${ f}\in {\mathfrak{D}}_{{{ p}}_{{ n}},{{ q}}_{{ n}}}\left({{ \zeta}}_{{ 0}}\right)$ for all ${ n}\in {\mathbb N}$,

\item  ${{ sup}}_{{ z}\in { K}}\left|{\left[{{{ p}}_{{ n}}}/{{{ q}}_{{ n}}}\right]}_{{ f},{{ \zeta}}_{{ 0}}}\left({ z}\right){ -}P\left({ z}\right)\right|\to { 0}$, as ${ n}\to { +}\infty $ and

\item  for every compact set$\ \ { L}'\subset \overline{\Omega}$ we have
\[\;\;{{ sup}}_{{\rm \ }{ z}\in { L}'}\left|{\left[{{{ p}}_{{ n}}}/{{{ q}}_{{ n}}}\right]}_{{ f},{{ \zeta}}_{{ 0}}}\left({ z}\right){ -}{ f}\left({ z}\right)\right|\to { 0},\;\; \text{as ${ n}\to { +}\infty $.}\]
\end{enumerate}
The set of all such functions ${ f}\in { A}\left({ \Omega}\right)$ is dense and ${{ G}}_{{ \delta}}$ in$\ \ { A}\left({ \Omega}\right)$.
\end{thm}
Let ${ K}$ vary in$\ \ {{ K}}_{{ m}}$,$\ \ { m}{ =}{ 1},{ 2}{ ,\dots }$ according to Lemma 4.2, ${ L}{ =}{{ L}}_{{ m}}$, ${ m}{ =}{ 1},{ 2}{ ,\dots }$ to be an exhausting family of compact subsets of ${ \Omega}$ and ${ L}'\equiv {{ L}}^{{ '}}_{{ m}}{ =}\left\{{ z}{ :\ }\left|{ z}\right|\le { m}\right\}\bigcap \overline{{ \Omega}}$. Then, applying again Baire's category theorem we obtain the
following.
\begin{thm}\label{thm4.4}
Let$\ \ \mathfrak{F}$ be a subset of ${\mathbb N}\times {\mathbb N}$ containing a sequence$\ \ \left({\widetilde{{ p}}}_{{ n}},{\widetilde{{ q}}}_{{ n}}\right)\in \mathfrak{F}$,$\ \ { n}{ =}{ 1},{ 2}{ ,\dots }$ with ${\widetilde{{ p}}}_{{ n}}\to { +}\infty $. Let ${ \Omega}\subset {\mathbb C}$ be a domain, such that $\left\{\infty \right\}\bigcup \left({\mathbb C}\backslash \overline{{ \Omega}}\right)$ is connected and$\ \ {\overline{{ \Omega}}}^{{ o}}{ =}{ \Omega}$. Then there exists a holomorphic function$\ \ { f}\in { A}\left({ \Omega}\right)$, satisfying the following.

For every compact set ${ K}\subset {\mathbb C}\backslash \overline{{ \Omega}}$ with connected complement and every polynomial ${ P}$, there exists a sequence$\ \ \left({{ p}}_{{ n}},{{ q}}_{{ n}}\right)\in \mathfrak{F}$,$\ { n}{ =}{ 1},{ 2}{ ,\dots }$,
    such that for every compact set ${L}$ contained in ${\Omega}$ we have:
\begin{enumerate}
\item  ${ f}\in {\mathfrak{D}}_{{{ p}}_{{ n}},{{ q}}_{{ n}}}\left({ \zeta}\right)$ for all ${ n}\in {\mathbb N}$ and ${ \zeta}\in { L}$

\item ${{ sup}}_{{ \zeta}\in { L}{ ,\ }{ z}\in { K}}\left|{\left[{{{ p}}_{{ n}}}/{{{ q}}_{{ n}}}\right]}_{{ f},{ \zeta}}\left({ z}\right){ -}{ P}\left({ z}\right)\right|\to { 0}$, as ${ n}\to { +}\infty $.

\item  For every compact set$\ \ { L}'\subset \overline{{ \Omega}}$, we have
\[{{ sup}}_{{\rm  }{ \zeta}\in { L}{ ,\  }{ z}\in { L}'}\left|{\left[{{{ p}}_{{ n}}}/{{{ q}}_{{ n}}}\right]}_{{ f},{ \zeta}}\left({ z}\right){ -}{ f}\left({ z}\right)\right|\to { 0},\; \text{as}\; \;{ n}\to { +}\infty.\]
\end{enumerate}
The set of all such functions ${ f}\in { A}\left({ \Omega}\right)$ is dense and ${{ G}}_{{ \delta}}$ in$\ \ { A}\left({ \Omega}\right)$.
\end{thm}
Notice that the complement in ${\mathbb C}\bigcup \left\{\infty \right\}$ of the compact set ${{ L}}^{{ '}}_{{ m}}{ =}\left\{{ z}{ :\ }\left|{ z}\right|\le { m}\right\}\bigcap \overline{{ \Omega}}$ is equal to $\left[\left\{{ z}{ :\ }\left|{ z}\right|{ >}m\right\}\bigcup \left\{\infty \right\}\right]\bigcup \left[\left\{\infty \right\}\bigcup \left({\mathbb C}\backslash \overline{{ \Omega}}\right)\right]$ which is connected, as the union of two connected sets with common point $ \infty $. The interior of ${{ L}}^{{ '}}_{{ m}}$ is$\ \ { \Omega}\bigcap \left\{{ z}{ :\ }\left|{ z}\right|{ <}m\right\}$, since$\ \ {\overline{{ \Omega}}}^{{ o}}{ =}{ \Omega}$.\textbf{}
\begin{rem}\label{rem4.5}
In Theorems 4.3 and 4.4, we can obtain ${\left[{{{ p}}_{{ n}}}/{{{ q}}_{{ n}}}\right]}_{{ f},{ \zeta}}\left({ z}\right)\to { P}\left({ z}\right)$ compactly on every simply connected domain ${ G}\subset {\mathbb C}\backslash \overline{{ \Omega}}$ (\cite{Cos}). Thus, the approximation is valid for all derivatives as well, according to Weierstrass theorem. However, this is not possible anymore if we allow that ${P}$ is not a polynomial, but ${P}\in A({K})$ (see Remark \ref{rem3.6}).
\end{rem}

\noindent
{\bf Acknowledgement}: We acknowledge a helpful communication with Franck Wielonsky.

\end{document}